\documentclass[twoside,11pt]{article}
\usepackage{fullpage}

\usepackage{epsf,amsmath,graphics}
\usepackage[psamsfonts]{amssymb,eucal}

\def\temp{1.35}%
\let\tempp=\relax
\expandafter\ifx\csname psboxversion\endcsname\relax
\else
    \ifdim\temp cm>\psboxversion cm
    \else
      \let\temp=\psboxversion
      \let\tempp= 
    \fi
\fi
\tempp
\let\psboxversion=\temp
\catcode`\@=11
\def\psfortextures{
\def\PSspeci@l##1##2{%
\special{illustration ##1\space scaled ##2}%
}}%
\def\psfordvitops{
\def\PSspeci@l##1##2{%
\special{dvitops: import ##1\space \the\drawingwd \the\drawinght}%
}}%
\def\psfordvips{
\def\PSspeci@l##1##2{%
\d@my=0.1bp \d@mx=\drawingwd \divide\d@mx by\d@my
\includegraphics{##1\space}}}%
\def\psforoztex{
\def\PSspeci@l##1##2{%
\special{##1 \space
      ##2 1000 div dup scale
      \number-\psllx\space\space \number-\pslly\space\space translate
}}}%
\def\psfordvitps{
\def\dvitpsLiter@ldim##1{\dimen0=##1\relax
\special{dvitps: Literal "\number\dimen0\space"}}%
\def\PSspeci@l##1##2{%
\at(0bp;\drawinght){%
\special{dvitps: Include0 "psfig.psr"}
\dvitpsLiter@ldim{\drawingwd}%
\dvitpsLiter@ldim{\drawinght}%
\dvitpsLiter@ldim{\psllx bp}%
\dvitpsLiter@ldim{\pslly bp}%
\dvitpsLiter@ldim{\psurx bp}%
\dvitpsLiter@ldim{\psury bp}%
\special{dvitps: Literal "startTexFig"}%
\special{dvitps: Include1 "##1"}%
\special{dvitps: Literal "endTexFig"}%
}}}%
\def\psfordvialw{
\def\PSspeci@l##1##2{
\special{language "PostScript",
position = "bottom left",
literal "  \psllx\space \pslly\space translate
  ##2 1000 div dup scale
  -\psllx\space -\pslly\space translate",
include "##1"}
}}%
\def\psforptips{
\def\PSspeci@l##1##2{{
\d@mx=\psurx bp
\advance \d@mx by -\psllx bp
\divide \d@mx by 1000\multiply\d@mx by \xscale
\incm{\d@mx}
\let\tmpx\dimincm
\d@my=\psury bp
\advance \d@my by -\pslly bp
\divide \d@my by 1000\multiply\d@my by \xscale
\incm{\d@my}
\let\tmpy\dimincm
\d@mx=-\psllx bp
\divide \d@mx by 1000\multiply\d@mx by \xscale
\d@my=-\pslly bp
\divide \d@my by 1000\multiply\d@my by \xscale
\at(\d@mx;\d@my){\special{ps:##1 x=\tmpx cm, y=\tmpy cm}}
}}}%
\def\psonlyboxes{
\def\PSspeci@l##1##2{%
\at(0cm;0cm){\boxit{\vbox to\drawinght
  {\vss\hbox to\drawingwd{\at(0cm;0cm){\hbox{({\tt##1})}}\hss}}}}
}}%
\def\psloc@lerr#1{%
\let\savedPSspeci@l=\PSspeci@l%
\def\PSspeci@l##1##2{%
\at(0cm;0cm){\boxit{\vbox to\drawinght
  {\vss\hbox to\drawingwd{\at(0cm;0cm){\hbox{({\tt##1}) #1}}\hss}}}}
\let\PSspeci@l=\savedPSspeci@l
}}%
\newread\pst@mpin
\newdimen\drawinght\newdimen\drawingwd
\newdimen\psxoffset\newdimen\psyoffset
\newbox\drawingBox
\newcount\xscale \newcount\yscale \newdimen\pscm\pscm=1cm
\newdimen\d@mx \newdimen\d@my
\newdimen\pswdincr \newdimen\pshtincr
\let\ps@nnotation=\relax
{\catcode`\|=0 |catcode`|\=12 |catcode`|
|catcode`#=12 |catcode`*=14
|xdef|backslashother{\}*
|xdef|percentother{
|xdef|tildeother{~}*
|xdef|sharpother{#}*
}%
\def\R@moveMeaningHeader#1:->{}%
\def\uncatcode#1{%
\edef#1{\expandafter\R@moveMeaningHeader\meaning#1}}%
\def\execute#1{#1}
\def\psm@keother#1{\catcode`#112\relax}
\def\executeinspecs#1{%
\execute{\begingroup\let\do\psm@keother\dospecials\catcode`\^^M=9#1\endgroup}}%
\def\@mpty{}%
\def\matchexpin#1#2{
  \fi%
  \edef\tmpb{{#2}}%
  \expandafter\makem@tchtmp\tmpb%
  \edef\tmpa{#1}\edef\tmpb{#2}%
  \expandafter\expandafter\expandafter\m@tchtmp\expandafter\tmpa\tmpb\endm@tch%
  \if\match%
}%
\def\matchin#1#2{%
  \fi%
  \makem@tchtmp{#2}%
  \m@tchtmp#1#2\endm@tch%
  \if\match%
}%
\def\makem@tchtmp#1{\def\m@tchtmp##1#1##2\endm@tch{%
  \def\tmpa{##1}\def\tmpb{##2}\let\m@tchtmp=\relax%
  \ifx\tmpb\@mpty\def\match{YN}%
  \else\def\match{YY}\fi%
}}%
\def\incm#1{{\psxoffset=1cm\d@my=#1
 \d@mx=\d@my
  \divide\d@mx by \psxoffset
  \xdef\dimincm{\number\d@mx.}
  \advance\d@my by -\number\d@mx cm
  \multiply\d@my by 100
 \d@mx=\d@my
  \divide\d@mx by \psxoffset
  \edef\dimincm{\dimincm\number\d@mx}
  \advance\d@my by -\number\d@mx cm
  \multiply\d@my by 100
 \d@mx=\d@my
  \divide\d@mx by \psxoffset
  \xdef\dimincm{\dimincm\number\d@mx}
}}%
\newif\ifNotB@undingBox
\newhelp\PShelp{Proceed: you'll have a 5cm square blank box instead of
your graphics.}%
\def\s@tsize#1 #2 #3 #4\@ndsize{
  \def\psllx{#1}\def\pslly{#2}%
  \def\psurx{#3}\def\psury{#4}
  \ifx\psurx\@mpty\NotB@undingBoxtrue
  \else
    \drawinght=#4bp\advance\drawinght by-#2bp
    \drawingwd=#3bp\advance\drawingwd by-#1bp
  \fi
  }%
\def\sc@nBBline#1:#2\@ndBBline{\edef\p@rameter{#1}\edef\v@lue{#2}}%
\def\g@bblefirstblank#1#2:{\ifx#1 \else#1\fi#2}%
{\catcode`\%=12
\xdef\B@undingBox{
\def\ReadPSize#1{
 \readfilename#1\relax
 \let\PSfilename=\lastreadfilename
 \openin\pst@mpin=#1\relax
 \ifeof\pst@mpin \errhelp=\PShelp
   \errmessage{I haven't found your postscript file (\PSfilename)}%
   \psloc@lerr{was not found}%
   \s@tsize 0 0 142 142\@ndsize
   \closein\pst@mpin
 \else
   \if\matchexpin{\GlobalInputList}{, \lastreadfilename}%
   \else\xdef\GlobalInputList{\GlobalInputList, \lastreadfilename}%
     \immediate\write\psbj@inaux{\lastreadfilename,}%
   \fi%
   \loop
     \executeinspecs{\catcode`\ =10\global\read\pst@mpin to\n@xtline}%
     \ifeof\pst@mpin
       \errhelp=\PShelp
       \errmessage{(\PSfilename) is not an Encapsulated PostScript File:
           I could not find any \B@undingBox: line.}%
       \edef\v@lue{0 0 142 142:}%
       \psloc@lerr{is not an EPSFile}%
       \NotB@undingBoxfalse
     \else
       \expandafter\sc@nBBline\n@xtline:\@ndBBline
       \ifx\p@rameter\B@undingBox\NotB@undingBoxfalse
         \edef\t@mp{%
           \expandafter\g@bblefirstblank\v@lue\space\space\space}%
         \expandafter\s@tsize\t@mp\@ndsize
       \else\NotB@undingBoxtrue
       \fi
     \fi
   \ifNotB@undingBox\repeat
   \closein\pst@mpin
 \fi
\message{#1}%
}%
\def\psboxto(#1;#2)#3{\vbox{%
   \ReadPSize{#3}%
   \advance\pswdincr by \drawingwd
   \advance\pshtincr by \drawinght
   \divide\pswdincr by 1000
   \divide\pshtincr by 1000
   \d@mx=#1
   \ifdim\d@mx=0pt\xscale=1000
         \else \xscale=\d@mx \divide \xscale by \pswdincr\fi
   \d@my=#2
   \ifdim\d@my=0pt\yscale=1000
         \else \yscale=\d@my \divide \yscale by \pshtincr\fi
   \ifnum\yscale=1000
         \else\ifnum\xscale=1000\xscale=\yscale
                    \else\ifnum\yscale<\xscale\xscale=\yscale\fi
              \fi
   \fi
   \divide\drawingwd by1000 \multiply\drawingwd by\xscale
   \divide\drawinght by1000 \multiply\drawinght by\xscale
   \divide\psxoffset by1000 \multiply\psxoffset by\xscale
   \divide\psyoffset by1000 \multiply\psyoffset by\xscale
   \global\divide\pscm by 1000
   \global\multiply\pscm by\xscale
   \multiply\pswdincr by\xscale \multiply\pshtincr by\xscale
   \ifdim\d@mx=0pt\d@mx=\pswdincr\fi
   \ifdim\d@my=0pt\d@my=\pshtincr\fi
   \message{scaled \the\xscale}%
 \hbox to\d@mx{\hss\vbox to\d@my{\vss
   \global\setbox\drawingBox=\hbox to 0pt{\kern\psxoffset\vbox to 0pt{%
      \kern-\psyoffset
      \PSspeci@l{\PSfilename}{\the\xscale}%
      \vss}\hss\ps@nnotation}%
   \global\wd\drawingBox=\the\pswdincr
   \global\ht\drawingBox=\the\pshtincr
   \global\drawingwd=\pswdincr
   \global\drawinght=\pshtincr
   \baselineskip=0pt
   \copy\drawingBox
 \vss}\hss}%
  \global\psxoffset=0pt
  \global\psyoffset=0pt
  \global\pswdincr=0pt
  \global\pshtincr=0pt 
  \global\pscm=1cm 
}}%
\def\psboxscaled#1#2{\vbox{%
  \ReadPSize{#2}%
  \xscale=#1
  \message{scaled \the\xscale}%
  \divide\pswdincr by 1000 \multiply\pswdincr by \xscale
  \divide\pshtincr by 1000 \multiply\pshtincr by \xscale
  \divide\psxoffset by1000 \multiply\psxoffset by\xscale
  \divide\psyoffset by1000 \multiply\psyoffset by\xscale
  \divide\drawingwd by1000 \multiply\drawingwd by\xscale
  \divide\drawinght by1000 \multiply\drawinght by\xscale
  \global\divide\pscm by 1000
  \global\multiply\pscm by\xscale
  \global\setbox\drawingBox=\hbox to 0pt{\kern\psxoffset\vbox to 0pt{%
     \kern-\psyoffset
     \PSspeci@l{\PSfilename}{\the\xscale}%
     \vss}\hss\ps@nnotation}%
  \advance\pswdincr by \drawingwd
  \advance\pshtincr by \drawinght
  \global\wd\drawingBox=\the\pswdincr
  \global\ht\drawingBox=\the\pshtincr
  \global\drawingwd=\pswdincr
  \global\drawinght=\pshtincr
  \baselineskip=0pt
  \copy\drawingBox
  \global\psxoffset=0pt
  \global\psyoffset=0pt
  \global\pswdincr=0pt
  \global\pshtincr=0pt 
  \global\pscm=1cm
}}%
\def\psbox#1{\psboxscaled{1000}{#1}}%
\newif\ifn@teof\n@teoftrue
\newif\ifc@ntrolline
\newif\ifmatch
\newread\j@insplitin
\newwrite\j@insplitout
\newwrite\psbj@inaux
\immediate\openout\psbj@inaux=psbjoin.aux
\immediate\write\psbj@inaux{\string\joinfiles}%
\immediate\write\psbj@inaux{\jobname,}%
\def\toother#1{\ifcat\relax#1\else\expandafter%
  \toother@ux\meaning#1\endtoother@ux\fi}%
\def\toother@ux#1 #2#3\endtoother@ux{\def\tmp{#3}%
  \ifx\tmp\@mpty\def\tmp{#2}\let\next=\relax%
  \else\def\next{\toother@ux#2#3\endtoother@ux}\fi%
\next}%
\let\readfilenamehook=\relax
\def\re@d{\expandafter\re@daux}
\def\re@daux{\futurelet\nextchar\stopre@dtest}%
\def\re@dnext{\xdef\lastreadfilename{\lastreadfilename\nextchar}%
  \afterassignment\re@d\let\nextchar}%
\def\stopre@d{\egroup\readfilenamehook}%
\def\stopre@dtest{%
  \ifcat\nextchar\relax\let\nextread\stopre@d
  \else
    \ifcat\nextchar\space\def\nextread{%
      \afterassignment\stopre@d\chardef\nextchar=`}%
    \else\let\nextread=\re@dnext
      \toother\nextchar
      \edef\nextchar{\tmp}%
    \fi
  \fi\nextread}%
\def\readfilename{\bgroup%
  \let\\=\backslashother \let\%=\percentother \let\~=\tildeother
  \let\#=\sharpother \xdef\lastreadfilename{}%
  \re@d}%
\xdef\GlobalInputList{\jobname}%
\def\psnewinput{%
  \def\readfilenamehook{
    \if\matchexpin{\GlobalInputList}{, \lastreadfilename}%
    \else\xdef\GlobalInputList{\GlobalInputList, \lastreadfilename}%
      \immediate\write\psbj@inaux{\lastreadfilename,}%
    \fi%
    \let\readfilenamehook=\relax%
    \ps@ldinput\lastreadfilename\relax%
  }\readfilename%
}%
\expandafter\ifx\csname @@input\endcsname\relax    
  \immediate\let\ps@ldinput=\input\def\input{\psnewinput}%
\else
  \immediate\let\ps@ldinput=\@@input
  \def\@@input{\psnewinput}%
\fi%
\def\nowarnopenout{%
 \def\warnopenout##1##2{%
   \readfilename##2\relax
   \message{\lastreadfilename}%
   \immediate\openout##1=\lastreadfilename\relax}}%
\def\warnopenout#1#2{%
 \readfilename#2\relax
 \def\t@mp{TrashMe,psbjoin.aux,psbjoint.tex,}\uncatcode\t@mp
 \if\matchexpin{\t@mp}{\lastreadfilename,}%
 \else
   \immediate\openin\pst@mpin=\lastreadfilename\relax
   \ifeof\pst@mpin
     \else
     \edef\tmp{{If the content of this file is precious to you, this
is your last chance to abort (ie press x or e) and rename it before
retexing (\jobname). If you're sure there's no file
(\lastreadfilename) in the directory of (\jobname), then go on: I'm
simply worried because you have another (\lastreadfilename) in some
directory I'm looking in for inputs...}}%
     \errhelp=\tmp
     \errmessage{I may be about to replace your file named \lastreadfilename}%
   \fi
   \immediate\closein\pst@mpin
 \fi
 \message{\lastreadfilename}%
 \immediate\openout#1=\lastreadfilename\relax}%
{\catcode`\%=12\catcode`\*=14
\gdef\splitfile#1{*
 \readfilename#1\relax
 \immediate\openin\j@insplitin=\lastreadfilename\relax
 \ifeof\j@insplitin
   \message{! I couldn't find and split \lastreadfilename!}*
 \else
   \immediate\openout\j@insplitout=TrashMe
   \message{< Splitting \lastreadfilename\space into}*
   \loop
     \ifeof\j@insplitin
       \immediate\closein\j@insplitin\n@teoffalse
     \else
       \n@teoftrue
       \executeinspecs{\global\read\j@insplitin to\spl@tinline\expandafter
         \ch@ckbeginnewfile\spl@tinline
       \ifc@ntrolline
       \else
         \toks0=\expandafter{\spl@tinline}*
         \immediate\write\j@insplitout{\the\toks0}*
       \fi
     \fi
   \ifn@teof\repeat
   \immediate\closeout\j@insplitout
 \fi\message{>}*
}*
\gdef\ch@ckbeginnewfile#1
 \def\t@mp{#1}*
 \ifx\@mpty\t@mp
   \def\t@mp{#3}*
   \ifx\@mpty\t@mp
     \global\c@ntrollinefalse
   \else
     \immediate\closeout\j@insplitout
     \warnopenout\j@insplitout{#2}*
     \global\c@ntrollinetrue
   \fi
 \else
   \global\c@ntrollinefalse
 \fi}*
\gdef\joinfiles#1\into#2{*
 \message{< Joining following files into}*
 \warnopenout\j@insplitout{#2}*
 \message{:}*
 {*
 \edef\w@##1{\immediate\write\j@insplitout{##1}}*
\w@{
\w@{
\w@{
\w@{
\w@{
\w@{
\w@{
\w@{
\w@{
\w@{
\w@{\string\input\space psbox.tex}*
\w@{\string\splitfile{\string\jobname}}*
\w@{\string\let\string\autojoin=\string\relax}*
}*
 \expandafter\tre@tfilelist#1, \endtre@t
 \immediate\closeout\j@insplitout
 \message{>}*
}*
\gdef\tre@tfilelist#1, #2\endtre@t{*
 \readfilename#1\relax
 \ifx\@mpty\lastreadfilename
 \else
   \immediate\openin\j@insplitin=\lastreadfilename\relax
   \ifeof\j@insplitin
     \errmessage{I couldn't find file \lastreadfilename}*
   \else
     \message{\lastreadfilename}*
     \immediate\write\j@insplitout{
     \executeinspecs{\global\read\j@insplitin to\oldj@ininline}*
     \loop
       \ifeof\j@insplitin\immediate\closein\j@insplitin\n@teoffalse
       \else\n@teoftrue
         \executeinspecs{\global\read\j@insplitin to\j@ininline}*
         \toks0=\expandafter{\oldj@ininline}*
         \let\oldj@ininline=\j@ininline
         \immediate\write\j@insplitout{\the\toks0}*
       \fi
     \ifn@teof
     \repeat
   \immediate\closein\j@insplitin
   \fi
   \tre@tfilelist#2, \endtre@t
 \fi}*
}%
\def\autojoin{%
 \immediate\write\psbj@inaux{\string\into{psbjoint.tex}}%
 \immediate\closeout\psbj@inaux
 \expandafter\joinfiles\GlobalInputList\into{psbjoint.tex}%
}%
\def\centinsert#1{\midinsert\line{\hss#1\hss}\endinsert}%
\def\psannotate#1#2{\vbox{%
  \def\ps@nnotation{#2\global\let\ps@nnotation=\relax}#1}}%
\def\pscaption#1#2{\vbox{%
   \setbox\drawingBox=#1
   \copy\drawingBox
   \vskip\baselineskip
   \vbox{\hsize=\wd\drawingBox\setbox0=\hbox{#2}%
     \ifdim\wd0>\hsize
       \noindent\unhbox0\tolerance=5000
    \else\centerline{\box0}%
    \fi
}}}%
\def\at(#1;#2)#3{\setbox0=\hbox{#3}\ht0=0pt\dp0=0pt
  \rlap{\kern#1\vbox to0pt{\kern-#2\box0\vss}}}%
\newdimen\gridht \newdimen\gridwd
\def\gridfill(#1;#2){%
  \setbox0=\hbox to 1\pscm
  {\vrule height1\pscm width.4pt\leaders\hrule\hfill}%
  \gridht=#1
  \divide\gridht by \ht0
  \multiply\gridht by \ht0
  \gridwd=#2
  \divide\gridwd by \wd0
  \multiply\gridwd by \wd0
  \advance \gridwd by \wd0
  \vbox to \gridht{\leaders\hbox to\gridwd{\leaders\box0\hfill}\vfill}}%
\def\fillinggrid{\at(0cm;0cm){\vbox{%
  \gridfill(\drawinght;\drawingwd)}}}%
\def\textleftof#1:{%
  \setbox1=#1
  \setbox0=\vbox\bgroup
    \advance\hsize by -\wd1 \advance\hsize by -2em}%
\def\textrightof#1:{%
  \setbox0=#1
  \setbox1=\vbox\bgroup
    \advance\hsize by -\wd0 \advance\hsize by -2em}%
\def\endtext{%
  \egroup
  \hbox to \hsize{\valign{\vfil##\vfil\cr%
\box0\cr%
\noalign{\hss}\box1\cr}}}%
\def\frameit#1#2#3{\hbox{\vrule width#1\vbox{%
  \hrule height#1\vskip#2\hbox{\hskip#2\vbox{#3}\hskip#2}%
        \vskip#2\hrule height#1}\vrule width#1}}%
\def\boxit#1{\frameit{0.4pt}{0pt}{#1}}%
\catcode`\@=12 
\psfordvips   

\def\IMSmarkvadjust{0 pt}
\def\IMSmarkhadjust{0 pt}
\def\IMSmarkhpadding{0 pt}
\def\SBIMSMark#1#2#3{
 \font\SBF=cmss10 at 10 true pt
 \font\SBI=cmssi10 at 10 true pt
 \setbox0=\hbox{\SBF \hbox to \IMSmarkhpadding{\relax}
                Stony Brook IMS Preprint \##1}
 \setbox2=\hbox to \wd0{\hfil \SBI #2}
 \setbox4=\hbox to \wd0{\hfil \SBI #3}
 \setbox6=\hbox to \wd0{\hss
             \vbox{\hsize=\wd0 \parskip=0pt \baselineskip=10 true pt
                   \copy0 \break%
                   \copy2 \break%
                   \copy4 \break}}
 \dimen0=\ht6   \advance\dimen0 by \vsize \advance\dimen0 by 8 true pt
                \advance\dimen0 by -\pagetotal
	        \advance\dimen0 by \IMSmarkvadjust
 \dimen2=\hsize \advance\dimen2 by .25 true in
	        \advance\dimen2 by \IMSmarkhadjust

  \openin2=publishd.tex
  \ifeof2\setbox0=\hbox to 0pt{}
  \else 
     \setbox0=\hbox to 3.1 true in{
                \vbox to \ht6{\hsize=3 true in \parskip=0pt  \noindent  
                {\SBI Published in modified form:}\hfil\break
                \input publishd.tex 
                \vfill}}
  \fi
  \closein2
  \ht0=0pt \dp0=0pt
 \ht6=0pt \dp6=0pt
 \setbox8=\vbox to \dimen0{\vfill \hbox to \dimen2{\copy0 \hss \copy6}}
 \ht8=0pt \dp8=0pt \wd8=0pt
 \copy8
 \message{*** Stony Brook IMS Preprint #1, #2. #3 ***}
}

\newtheorem{DEF}{\bf Definition}
\newtheorem{THEO}{\bf Theorem}
\newtheorem{LEMMA}{\bf Lemma}
\newtheorem{CORO}{\bf Corollary}

\def\C{{{\mathbb C}}}
\def\B{{{\mathbb B}}}
\def\T{{{\mathbb T}}}
\def\Z{{{\mathbb Z}}}
\def\D{{{\mathbb D}}}
\def\R{{{\mathbb R}}}
\def\N{{{\mathbb N}}}
\def\H{{{\mathbb H}}}
\def\M{{{\mathbb M}}}
\def\P{{{\mathbb P}}}
\def\P{{{\mathbb P}^1}}
\def\proof{\noindent{\sc Proof. }}
\def\qedprop{\hfill{\hbox{%
  \hskip 1pt%
  \vrule width 7pt height 6pt depth 1.5pt%
  \hskip 1pt}}}
\def\qedlem{\hfill{$\square$}}
\def\re{{\rm Re}}
\def\im{{\rm Im}}

\pagestyle{myheadings}
\markboth{\centerline{\sc Xavier  Buff and Christian 
Henriksen}}{\centerline{Scaling ratios and triangles in Siegel 
disks.}}

\begin{document}
\thispagestyle{empty}
\vskip 2cm

\SBIMSMark{1999/4}{April 1999}{}

\centerline{\Large Scaling ratios and triangles in Siegel disks.}
\centerline{
\begin{tabular}{ccc}
\rule{7cm}{0cm}  & by & \rule{7cm}{0cm} \\
{\sc \large Xavier Buff} & &{\sc \large Christian Henriksen.} \\
{\small Universit\'e Paul Sabatier} 
  & & {\small The Technical University of Denmark} \\
{\small Laboratoire Emile Picard} & and 
& {\small Department of Mathematics} \\
{\small 31062 Toulouse Cedex} & & {\small 2800 Lyngby} \\
{\small France} & & {\small Denmark} \\
\end{tabular}}

\vskip .5cm
Let $f(z)=e^{2i\pi\theta} z+z^2$, where $\theta$ is a quadratic 
irrational. McMullen proved that the Siegel disk for $f$ is 
self-similar about the critical point. We give a lower bound for the 
ratio of self-similarity, and we show that if $\theta=(\sqrt 5-1)/2$ 
is the golden mean, then there exists a triangle contained in the 
Siegel disk, and with one vertex at the critical point. This answers a 
15 year old conjecture.
\vskip 1cm

\section{Introduction}

\begin{DEF}
The polynomial $P_\theta$ is defined by 
$$P_\theta(z)=e^{2i\pi\theta}z+z^2,$$
where $\theta$ has continued fraction expansion
$$\theta=[a_1,a_2,a_3,\ldots]=
\cfrac{1}{a_1+\cfrac{1}{a_2+\cfrac{1}{a_3+\ddots}}}.$$
\end{DEF}

In the following, for $x\in \R/\Z$, $\{x\}$ denotes the unique real number 
representing $x$ in $]-1/2,1/2]$, and 
$$\frac{p_n}{q_n}=[a_1,\ldots,a_n]$$
denote the rational approximation to $\theta$ obtained by truncating 
its continued fraction.

In 1942, Siegel \cite{si} proved that when $\theta$ is a diophantine 
number, the polynomial $P_\theta$ is conformally 
conjugate to a rotation near the origin. The maximal domain $D$ on which 
this conjugacy is defined is called the Siegel disk for $P_\theta$. 
It is the Fatou component of $P_\theta$ containing $0$. In particular, 
this result holds when $\theta$ is of bounded type, 
i.e. $\sup a_i<\infty$.

In 1986, Herman \cite{h} and \'Swi\c{a}tek \cite{sw}
proved that when $\theta$ is of bounded type, the boundary 
$\partial D$ of the Siegel disk is a quasi-circle containing the 
critical point $\omega_\theta=-e^{2i\pi\theta}/2$. The proof is based on a
quasi-conformal surgery due to Ghys and Douady (see \cite{d}).
In 1993, Petersen 
\cite{p} proved that the Julia set $J(P_\theta)$
has Lebesgue measure zero, and is locally connected. 

In 1997, McMullen \cite{mcm} obtained results concerning the geometry of
the Julia set $J(P_\theta)$. In particular, he proved that 
if $\theta$ is a quadratic irrational, then the boundary of the Siegel 
disk for $P_\theta$ is self-similar about the critical point.
This result was conjectured and observed 
numerically more than a decade ago by Manton, Nauenberg and Widom 
\cite{mn} \cite{w}. 

The number $\theta$ is a quadratic irrational if and only if 
the continued fraction of $\theta$ is preperiodic. In that case,
the rotation $x\mapsto x+\theta$, $x\in \R/\Z$ is self-similar
(see \cite{mcm} theorem 2.1). More precisely, if
$\theta=[a_1,a_2,\ldots]$, where $a_{n+s}=a_n$ for $n\geq N$, we can
set 
$$\alpha=\theta_{N+1}\theta_{N+2}\ldots\theta_{N+s},$$
where 
$\theta_i=[a_i,a_{i+1},a_{i+2},\ldots]$. Then for $n\geq N$,
$$\{q_{n+s}\theta\}=(-1)^s \alpha\cdot \{q_n\theta\}.$$

Our first goal is to prove the following result. 
\begin{THEO}
Let $\theta=[a_1,a_2,\ldots]$, where $a_{n+s}=a_n$ for $n\geq N$, 
be a quadratic irrational and $\lambda\in\D^*$ be the 
scaling ratio for the self-similarity of the Siegel disk of $P_\theta$ 
about the critical point. Besides, let 
$\alpha$ be defined as above. Then
$$0<\alpha<|\lambda|<1.$$
\end{THEO}
McMullen mentioned to us that this bound on $\lambda$ in terms of $\alpha$ 
via a modulus estimate is very similar to Bers' inequality for
quasifuchsian groups; there one knows that the
length of a hyperbolic geodesic in $Q(X,Y)$ is bounded by
the hyperbolic length of the corresponding geodesic on X or Y (see 
\cite{b} Theorem 3 and \cite{mcm3} Prop. 6.4).

In \cite{mcm} (corollary 7.5), 
McMullen also shows that when the continued fraction 
expansion of $\theta$ has odd period, then the boundary of the Siegel 
disk does not spiral about the critical point. This means that any 
continuous branch of ${\rm arg}(z-\omega_\theta)$ defined along 
$\partial D\setminus \{\omega_\theta\}$ is bounded. In particular, 
this result holds for the {\it golden mean} Siegel disk, where 
$\theta=(\sqrt 5-1)/2=[1,1,1,\ldots]$. Our second result is the 
following.
\begin{THEO}
Using the same notations, if $-\pi/\log(\alpha^2)>1/2$, then 
the Siegel disk of the polynomial $P_\theta$, 
contains a triangle with one vertex at $\omega_\theta$.
\end{THEO}

\begin{CORO}
The Siegel disk of the polynomial $P_\theta$, $\theta=(\sqrt 5-1)/2$,
contains a triangle with one vertex at $\omega_\theta$.
\end{CORO}

\begin{figure}[htbp]
\centerline{
\boxit{\psboxscaled{500}{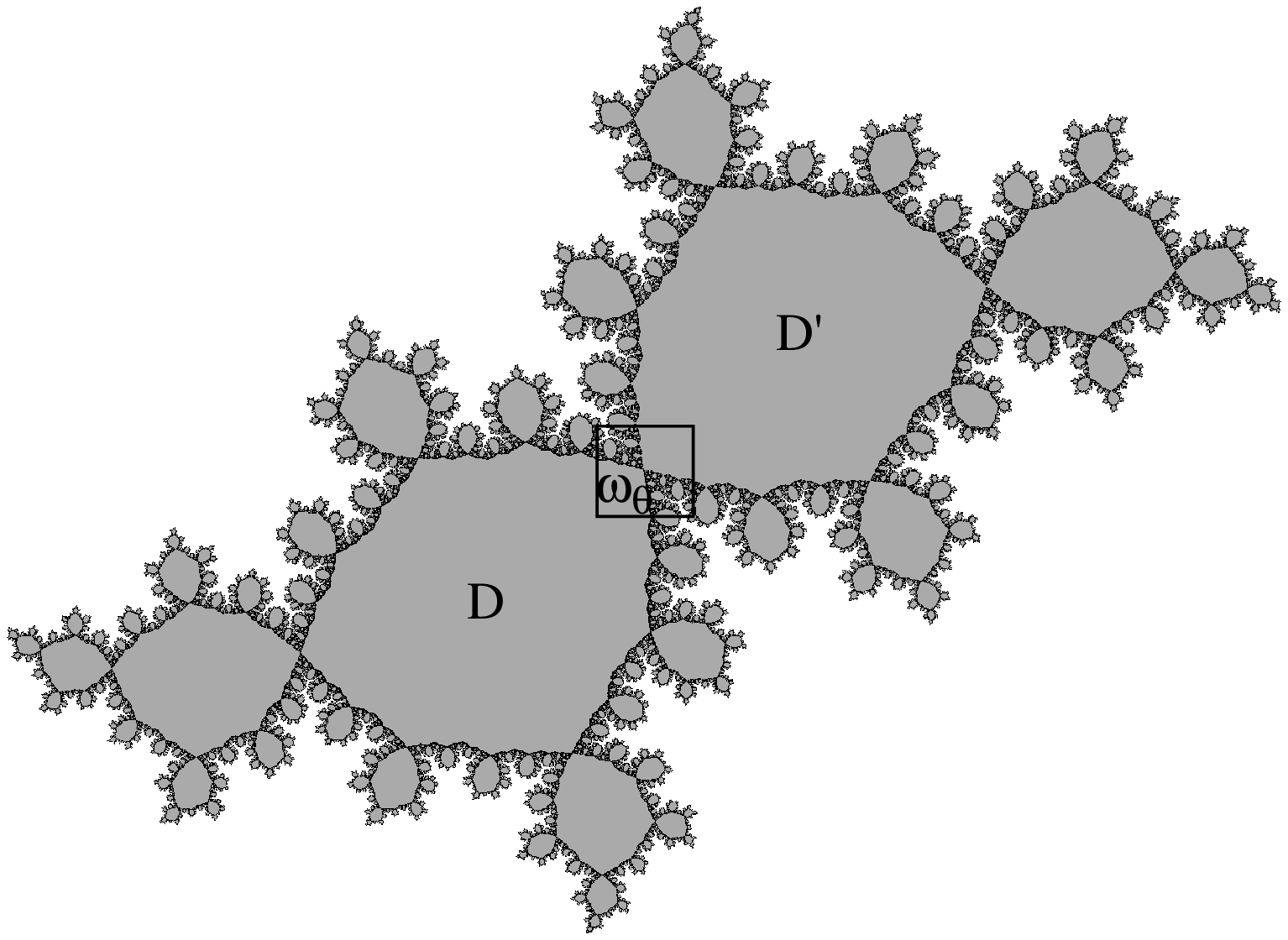}}
\boxit{\psboxscaled{500}{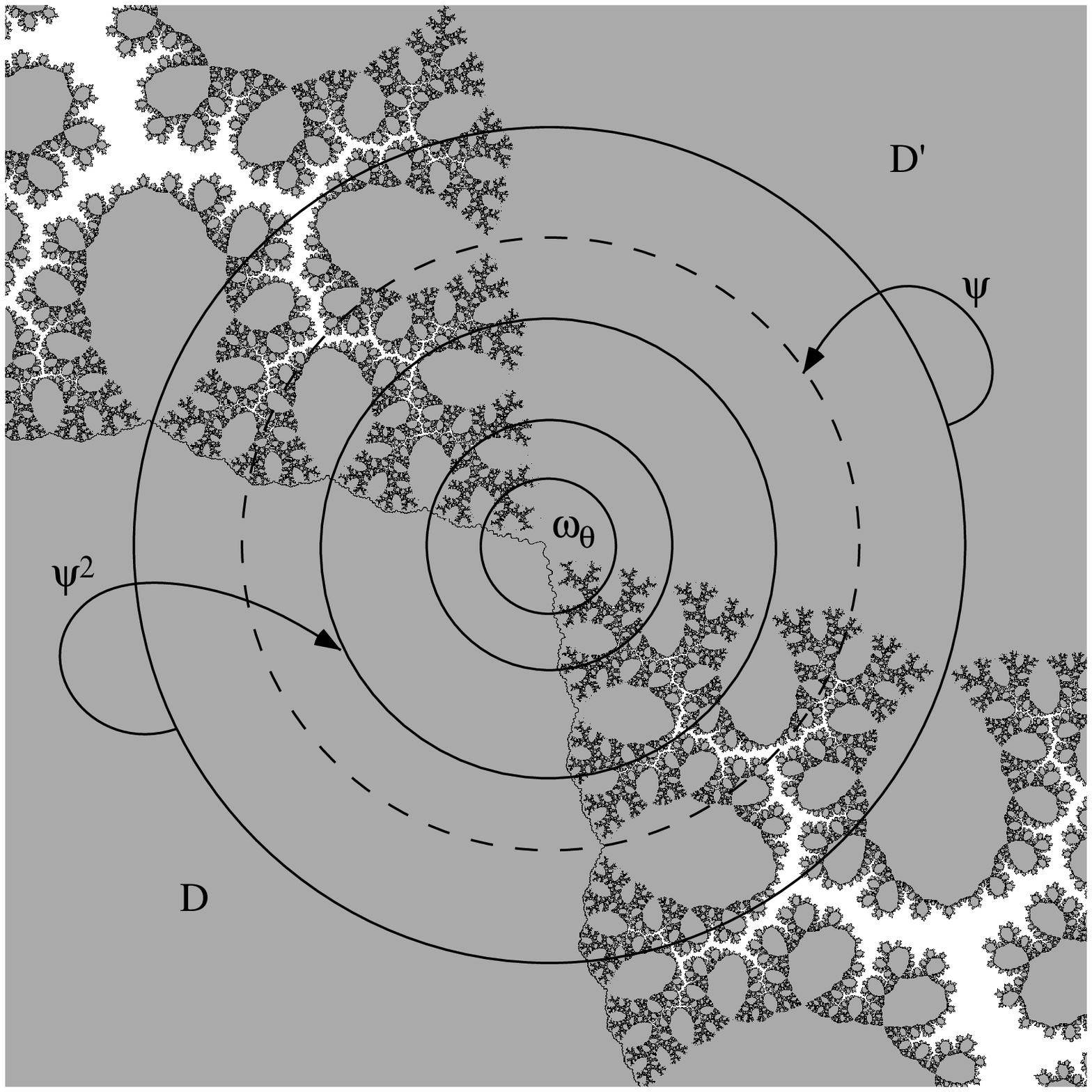}}}
\caption{The filled-in Julia set of the polynomial $P_\theta$, 
$\theta=(\sqrt 5-1)/2$.}
\label{figure1}
\end{figure}

The corollary is immediate since for $\theta=(\sqrt 
5-1)/2=[1,1,1,\ldots]$ we have $\alpha=\theta=(\sqrt 5-1)/2$, and
$$-\frac{\pi}{\log(\alpha^2)}\sim 3.264251306> \frac{1}{2}.$$
On figure \ref{figure1}, we have drawn the filled-in Julia set of the 
polynomial $P_\theta$, $\theta=(\sqrt 5-1)/2$. We have also zoomed 
near the critical point $\omega_{\theta}$ to show the self-similarity 
of the boundary of the Siegel disk.

\vskip.5cm
Acknowledgements. We wish to thank Adrien Douady for mentioning this
problem. We are grateful to Curt T. McMullen for carefully reading a first
version of this paper and suggesting several improvements.
We also want to thank John H. Hubbard, Misha Lyubich and Carsten L.
Petersen for valuable comments and the Departments of Mathematics at Cornell
University, at the Technical University of Denmark and at Universit\'e Paul
Sabatier in Toulouse for hospitality during the research
that went into this paper.

\section{The scaling ratio.}

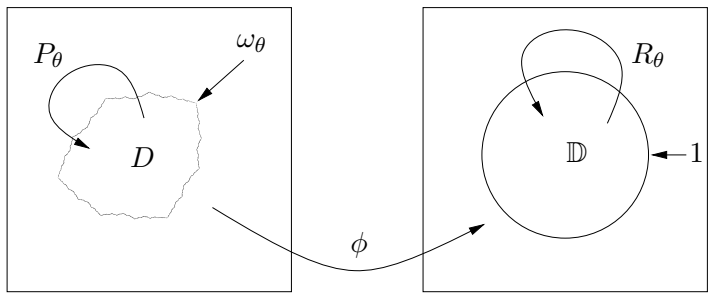
\begin{figure}[htbp]
\vskip.5cm
\centerline{
\begin{picture}(0,0)%
\psboxscaled{700}{Phi.pstex}%
\end{picture}%
\setlength{\unitlength}{0.000612500in}%
\begingroup\makeatletter\ifx\SetFigFont\undefined%
\gdef\SetFigFont#1#2#3#4#5{%
  \reset@font\fontsize{#1}{#2pt}%
  \fontfamily{#3}\fontseries{#4}\fontshape{#5}%
  \selectfont}%
\fi\endgroup%
\begin{picture}(6009,2457)(1114,-3991)
\put(1351,-1996){$P_\theta$}
\put(2181,-2896){$D$}
\put(3096,-1866){$\omega_\theta$}
\put(5926,-2851){$\D$}
\put(6981,-2851){$1$}
\put(4071,-3626){$\phi$}
\put(6481,-1996){$R_\theta$}
\end{picture}}
\caption{The linearizing map $\phi:\overline D\to \overline\D$ sending $\omega_\theta$ to
  $1$.}
\label{phi}
\end{figure}

In the following, $\theta=[a_1,a_2,\ldots]$, where $a_{n+s}=a_n$ for 
$n\geq N$, will always be a quadratic irrational. We denote by 
$R_\theta(z)=e^{2i\pi\theta} z$ the rotation of angle $\theta$. The polynomial 
$P_\theta$ has a Siegel disk $D$, and the conformal mapping 
$\phi:~D\to \D$, linearizes $P_\theta$, i.e. conjugates $P_\theta$ to 
the rotation $R_\theta$. By results of Herman and \'Swi\c{a}tek, 
this conjugacy extends to a quasi-symmetric conjugacy $\phi:~\partial 
D\to S^1$. Since $\phi$ is unique up to rotation, we can normalize it 
so that it maps $\omega_\theta\in\partial D$ to $1\in S^1$ (see Figure
\ref{phi}). 

Now, recall that
\begin{itemize}
\item[$\bullet$]{$p_n/q_n=[a_1,\ldots,a_n]$ is the rational approximation to 
$\theta$ obtained by truncating its continued fraction,}
\item[$\bullet$]{$\theta_i=[a_i,a_{i+1},a_{i+2},\ldots]$, and}
\item[$\bullet$]{$\alpha=\theta_{N+1}\theta_{N+2}\ldots\theta_{N+s}$.}
\end{itemize}
In \cite{mcm}, McMullen proves that for $n\geq N$,
$$\{q_{n+s}\theta\}=(-1)^s \alpha \{q_n\theta\}.$$
It follows that in a neighborhood of $z=1$, the contraction 
$$\begin{array}{l}
z\mapsto z^{\alpha}, \quad {\rm for}~s~{\rm even,~and}\\
z\mapsto \bar z^{\alpha}, \quad {\rm for}~s~{\rm odd}
\end{array}$$ 
conjugates $R_\theta^{q_n}$ to $R_\theta^{q_{n+s}}$, for $n\geq
N$. Let us prove it for $s$ odd. For $z$ in a sufficiently small
neighborhood of $1,$ we have
\begin{equation*}
\overline{R_\theta^{q_n}(z)}^\alpha  =  
  \overline{e^{2\pi i q_n \theta} z}^\alpha 
   =  \left( e^{-2 \pi i \{ q_n \theta \}} \bar{z} \right) ^\alpha 
   =  e^{ -2 \pi i \alpha \{ q_n \theta \} } \bar{z}^\alpha 
   =  R_\theta^{q_{n+s}}(\bar{z}^\alpha) .
\end{equation*}

In \cite{mcm} (theorem 7.1), McMullen proves that there is a 
neighborhood $U$ of $\omega_\theta$ and a constant $\varepsilon > 0$
such that for all $z\in U\cap \overline D$,
the  mapping $\psi$ defined by
$$\psi(z)=
\left\{\begin{array}{l}
\phi^{-1}\left([\phi(z)]^{\alpha}\right),\quad 
{\rm for}~s~{\rm even,~and}\\
\phi^{-1}\left([\overline{\phi(z)}]^{\alpha}\right),\quad 
{\rm for}~s~{\rm odd},
\end{array}
\right.$$
\begin{itemize}
\item[$\bullet$]{is well defined,}
\item[$\bullet$]{satisfies the expansion
$$\psi(z)=\left\{
\begin{array}{l}
\omega_\theta+\lambda(z-\omega_\theta)
   +{\cal O}\left(|z-\omega_{\theta}|^{1+\varepsilon}\right), \quad 
   {\rm for}~s~{\rm even,~or}\\
\omega_\theta+\lambda\overline{(z-\omega_\theta)}
   +{\cal O}\left(|z-\omega_{\theta}|^{1+\varepsilon}\right),\quad 
{\rm for}~s~{\rm odd},
\end{array}
\right.$$
with $0<|\lambda|<1$, and}
\item[$\bullet$]{conjugates $P_\theta^{q_n}$ to $P_\theta^{q_{n+s}}$.}
\end{itemize}
The main difficulty is to prove that 
$\psi$ is $C^{1+\varepsilon}$ at $\omega_\theta$.

Now let us define the scaling map
$$\Lambda(z)=\left\{
\begin{array}{l}
  \omega_\theta+\lambda(z-\omega_\theta),\quad {\rm for}~s~{\rm
  even,~or} \\
\omega_\theta+\lambda\overline{(z-\omega_\theta)},
\quad {\rm for}~s~{\rm odd},
\end{array}
\right.$$
For conveniency, we will use the spherical metric
$$d\sigma(z)=\frac{|dz|}{1+|z-\omega_\theta|^2}$$ 
(instead of the usual
$|dz|/(1+|z|^2)$). Then the distance between two points $x$ and $y$ in
$\P$ satisfy
$$\sigma(x,y) \leq \inf\left(|x-y|,\frac{1}{|x-\omega_\theta|}+
\frac{1}{|y-\omega_\theta|}\right).$$
This spherical metric enables us to define a Hausdorff distance $d_H$
between compact subsets of the sphere.
McMullen shows the following theorem.

\begin{THEO}{\sc McMullen (\cite{mcm}, corollary 7.3)}
The blow-ups $S_n=\Lambda^{-n}(\partial D)$ of the
boundary of the Siegel  
disk converge to a $\Lambda$-invariant quasi-circle through
$\infty$, for the  
Hausdorff topology on compact subsets of the sphere.
\end{THEO}

\proof Indeed, there exists constants $C_1$ and
$\delta>0$ such that for all $n$ large enough
$$d_H(S_n,S_{n+1})<C_1 |\lambda|^{n\delta}.$$
To see that, we need to prove that
\begin{itemize}
\item[$\bullet$]{for any $x\in S_n$ there exists a
$y$ in $S_{n+1}$ with $\sigma(x,y)< C_1 |\lambda|^{n\delta}$, and}
\item[$\bullet$]{for any $y$ in $S_{n+1}$ there exists a $x\in S_n$
with $\sigma(x,y)< C_1 |\lambda|^{n\delta}$.}
\end{itemize}

We will only prove the first point; a similar argument works for the
second one.
We first choose a constant $C$ and an open neighborhood
$U'$ of $\omega_\theta$ sufficiently small so that
$\psi(U')\subset U'$, and so that for any $z\in U'$,
$$|\psi(z)-\Lambda(z)| < C|z-\omega_\theta|^{1+\varepsilon}.$$
We then set 
$$\delta=\frac{1}{2}\left(1-\frac{1}{1+\varepsilon}\right),$$
which is positive. Observe that for all $n$ large enough, the
ball $B_n$ centered at $\omega_\theta$ with radius
$|\lambda|^{n(1-\delta)}$ is contained in $U'$. Then, for any $x\in
S_n$, $\Lambda^n(x)$ belongs to $\partial D$ and
\begin{itemize}
\item[$\bullet$]{either $\Lambda^n(x)\in B_n$, 
$\psi(\Lambda^n(x))\in \partial D$ and 
$y=\Lambda^{-(n+1)}\left(\psi(\Lambda^n(x)\right)$ belongs to
$S_{n+1}$; then a simple computation gives
$$\sigma(x,y)<\frac{C}{|\lambda|} |\lambda|^{n\delta},$$}
\item[$\bullet$]{or $\Lambda^n(x)\not\in B_n$ and
    $y=\Lambda^{-1}(x)$ belongs to $S_{n+1}$; moreover
$$\sigma(x,y)\leq\frac{1}{|x-\omega_\theta|}+\frac{1}{|y-\omega_\theta|}
\leq 2|\lambda|^{n\delta}.$$}
\end{itemize}
Hence $d_H(S_n, S_{n+1})$ is decreasing geometrically and the sequence
$S_n$ is converging for the Hausdorff topology to a limit ${\cal
  S}$ which has to be $\Lambda$-invariant. 
Since the sets $S_n$ are all $K$ quasi-circles with the same
$K$ (they are mapped onto each other by the scaling map $\Lambda$),
the limit is also a $K$ quasi-circle.
\qedprop

Since $P_\theta$ is a quadratic polynomial, the Siegel disk $D$ has 
one preimage $D'\neq D$ which is symmetric to $D$ with respect 
to $\omega_\theta$. The blow-ups $\Lambda^{-n}(D')$ and 
$\Lambda^{-n}(D)$ both converge, for the 
Hausdorff topology on compact subsets of the sphere, 
to $\Lambda$-invariant quasi-disks 
${\cal D}$ (bounded by the quasi-circle ${\cal S}$)
and ${\cal D}'$ passing through $\infty$ and $\omega_\theta$ (see figure 
\ref{figure1}). In particular, observe that ${\cal D}/\Lambda^2$ and
${\cal D}'/\Lambda^{2}$ are two  
annuli in the torus $(\C\setminus\{\omega_\theta\})/\Lambda^2$. 
We consider $\Lambda^2$ instead of $\Lambda$, because when $s$ is odd, 
$\Lambda$ is orientation reversing. Notice that when $s$ 
is even, this torus is conformally equivalent to $\C^*/\lambda^2$,
and when $s$ is odd, this torus is 
conformally equivalent to the torus $\C^*/(\lambda\bar \lambda)$.
Besides, the annuli are conformally equivalent. Let 
$$M={\rm mod}({\cal D}/\Lambda^2)={\rm mod}({\cal D}'/\Lambda^2)$$
be their modulus.

The key-point in this paper is that we can compute the exact value of
the modulus $M$. 
\begin{LEMMA}\label{modulus}
The modulus $M$ is equal to $-\pi/\log(\alpha^2)$.
\end{LEMMA}

\proof
Indeed, we can define the scaling map
$$A(z)=\left\{
\begin{array}{l}
1+ \alpha (z-1),\quad{\rm for}~s~{\rm even,~and}\\
1+ \alpha \overline{(z-1)},\quad {\rm for}~s~{\rm odd}.
\end{array}
\right.$$
It is the differential at $1$ of the contraction which conjugates 
$R_\theta^{q_n}$ to $R_\theta^{q_{n+s}}$, for $n\geq N$.
Then
$$\phi_n=A^{-n}\circ \phi\circ \Lambda^n$$
is a conformal equivalence between $\Lambda^{-n}(D)$ and 
$A^{-n}(\D)$, which extends quasi-symmetrically to a map 
$\phi_n:~\Lambda^{-n}(\partial D)\to A^{-n}(\partial\D)$. Besides,
$\phi_n(\omega_\theta)=1$, and
$$A\circ \phi_{n+1}=\phi_n\circ \Lambda.$$
By Caratheodory's convergence theorem, 
the sequence $\phi_n$ converges when $n$ tends to infinity, to a 
conformal map 
$$\phi_\infty:~{\cal D}\to \H=\{z\in \C~|~{\rm Re}(z)<1\}$$ 
such that $A\circ \phi_\infty=\phi_\infty\circ \Lambda$ (see 
\cite{mcm} Theorem 8.1, statement 7). In 
particular, we see that the annulus
${\cal D}/\Lambda^2$ is isomorphic to the annulus $\H/A^2$. This last 
annulus has a modulus $M=-\pi/\log(\alpha^2)$. 
\qedlem

We will now use a classical inequality on annuli 
embedded in a torus.

\begin{LEMMA}\label{severalAnnInTorus}
Let $A_i \subset \T$, be disjoint annuli embedded in a torus 
$$\T = \C/(2\pi i\Z + \tau \Z), \quad \mathrm{Re}(\tau) > 0.$$
The segment $[0,\tau]$ projects to a simple closed curve $\gamma$ on $\T$.
If the annuli $A_i$ are homotopic to $\gamma$, then
$$\sum_{i=1}^n {\rm mod}A_i \leq \frac{2\pi \re(\tau)}{|\tau|^2}.$$ 
\end{LEMMA}

\proof Let $B_i$ be the annulus 
$$\{z~|~0 < \im(z) < h_i\}/\Z,$$
where $\Z$ acts by translations, with $h_i = {\rm mod}B_i =
{\rm mod} A_i$, so that $A_i$ and $B_i$ are conformally equivalent. 
And let $f_i:~B_i \to A_i$ be a conformal mapping. We can 
endow the torus $\T$ with the Euclidean metric; then the simple closed curve
$$x \mapsto f_i(x,y), \qquad 0 \le x \le 1$$
has length at least $|\tau|$. Hence, we find
\begin{eqnarray*}
2\pi \re(\tau) &=& {\rm Area}(\T) \\
& \ge & \sum_i {\rm Area}(A_i) \\
& = & \sum_i \int_{B_i}|f_i'(x,y)|^2 dxdy \\
& = & \sum_i \int_0^{h_i} \left(\int_0^1 |f_i'(x,y)|^2 dx\right)dy \\
& \ge & \sum_i
\int_0^{h_i} \left(\int_0^1 |f_i'(x,y)| dx\right)^2 dy \\
&\ge& \sum_i \int_0^{h_i} |\tau|^2 dy \\
& = & \sum_i h_i |\tau|^2.
\end{eqnarray*}
This is the required inequality. 
\qedlem

In our case, $\tau$ is a branch of $\log(\lambda^2)$ when $s$ is even 
and of $\log(\lambda \bar \lambda)$ when $n$ is odd. Using 
$$|{\rm Re}(\tau)|=-\log(|\lambda|^2),$$
we get
$$2M\leq -\frac{2\pi}{\log(|\lambda|^2)}.$$
Combining this with the exact value of the modulus $M$ gives
$$-\frac{2\pi}{\log(\alpha^2)}\leq-\frac{2\pi}{\log(|\lambda|^2)},$$
which is equivalent to
$$\alpha<|\lambda|.$$
This proves theorem 1.

\section{Triangle in the golden Siegel disk.}

We will now show that if $\theta=(\sqrt 5-1)/2$ is the golden mean, then 
the Siegel disk of $P_\theta$ contains a triangle with vertex at the 
critical point $\omega_\theta$. It is enough to show that the quasi-disk 
${\cal D}$ contains a sector with vertex at $\omega_\theta$ . 

McMullen has already done the main step in 
that direction (see \cite{mcm} corollary 7.5). He proved that when $s$ is 
odd (and in our case $s=1$), the boundary of the Siegel disk does not 
spiral about the critical point. That means that there exists a continuous
branch of $\chi(z)=\log(z-\omega_\theta)$ defined on $\cal D$ with bounded 
imaginary part. Indeed, the condition that ${\cal D}$ is 
$\Lambda$-invariant implies that ${\cal D}$ is 
$\Lambda^2$-invariant. Now the scaling ratio of $\Lambda^2$ is
$\lambda \bar{\lambda}$ and consequently the strip $\chi({\cal D})$ is
invariant by 
the translation $T(z)=z+\log(\lambda \bar{\lambda})$. This translation being 
real, the imaginary part of $\chi(z)$ is bounded when $z\in {\cal D}$.

\begin{figure}[htbp]
\centerline{
\boxit{\psboxscaled{500}{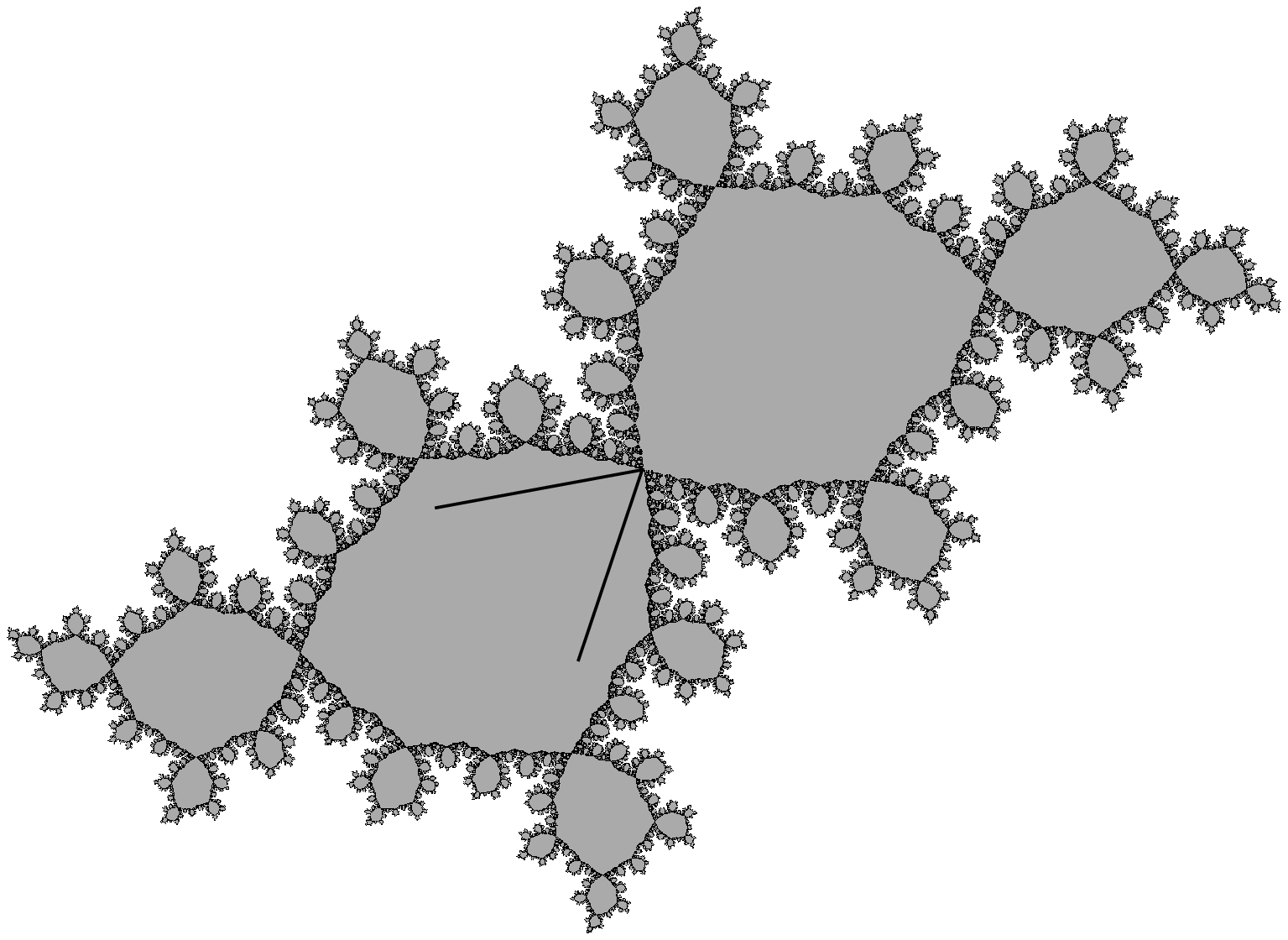}}
\boxit{\psboxscaled{500}{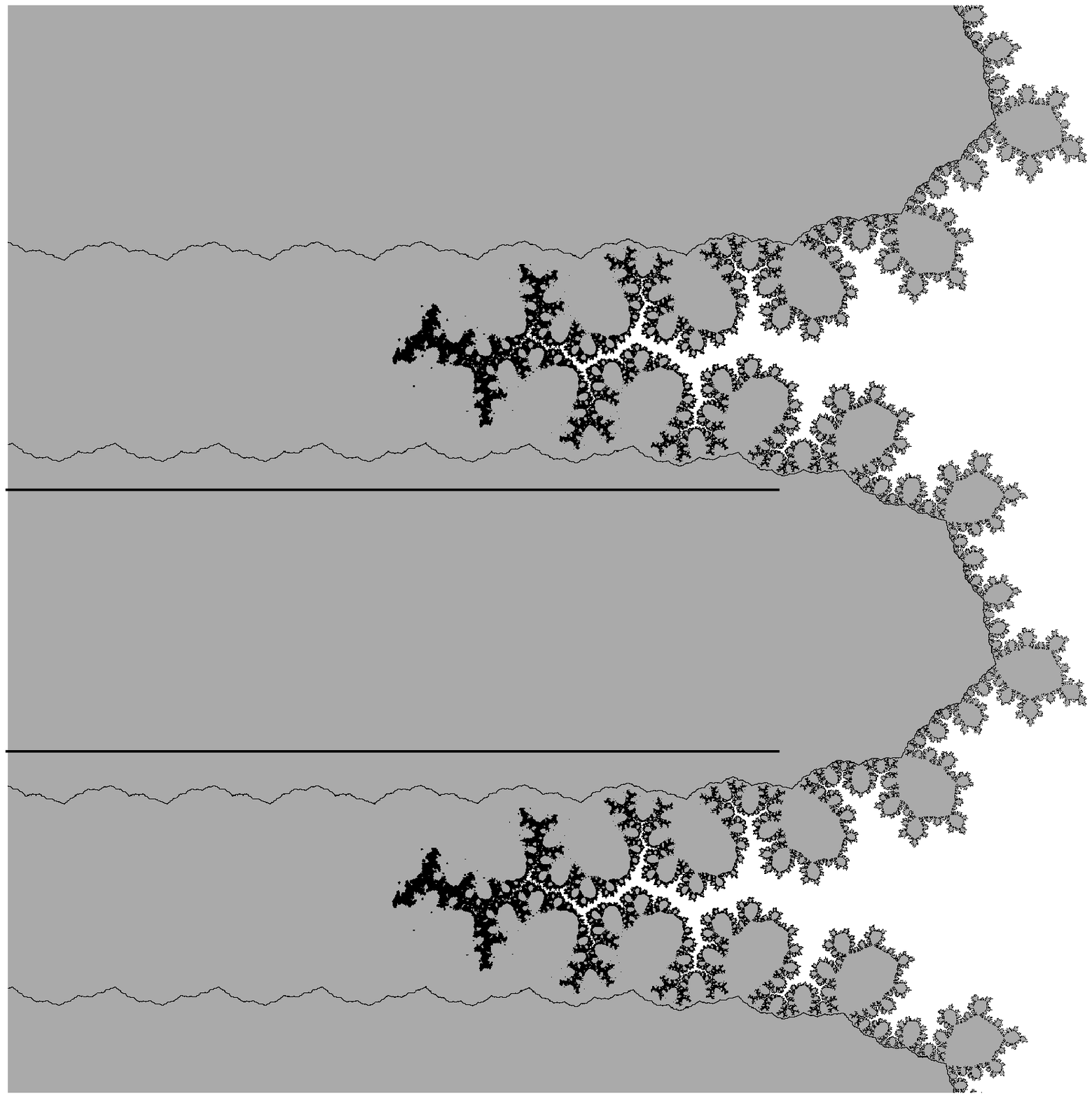}}}
\caption{The filled-in Julia set of the polynomial $P_\theta$, 
$\theta=(\sqrt 5-1)/2$ and its image under the map 
$\chi(z)=\log(z-\omega_{\theta})$.}
\label{figure3}
\end{figure}

To prove the existence of a sector in ${\cal D}$, it is enough to
show that the strip $\chi({\cal D})$ contains a 
horizontal band
$$B=\{z\in \C~|~y_1<{\rm Im}(z)<y_2\},$$
for some $y_1<y_2$ in $\R$ (see figure \ref{figure3}).

On figure \ref{figure3}, we have drawn the Julia set 
of the polynomial $P_\theta$, 
$\theta=(\sqrt 5-1)/2$ and its image under the map 
$\chi(z)=\log(z-\omega_{\theta})$. It is very difficult to get a good 
picture of the Julia set $J(P_\theta)$ 
near the critical point $\omega_\theta$. However, it is 
possible to get a good idea of the boundary of the Siegel disk since the 
orbit of the critical point is dense in $\partial D$.

To prove that the strip $\chi({\cal D})$ contains a 
horizontal band, recall that the quotient 
$\chi({\cal D})/T$ is an annulus of modulus $-\pi/\log(\alpha^2)$, where 
in our case $\alpha=\theta=(\sqrt 5-1)/2$. Since 
$$-\frac{\pi}{\log(\alpha^2)}\sim 3.264251306> \frac{1}{2} ,$$
the existence of a sector in ${\cal D}$ is a consequence of the 
following lemma (see \cite{mcm1}, Thm 2.1).

\begin{LEMMA}
Assume $S$ is a periodic strip of period $\tau\in \R$, i.e. $S+\tau=\tau$.
If $S/\tau$ is an annulus of modulus
$${\rm mod}(S/\tau)>\frac{1}{2},$$ then $S$ contains a horizontal band
$$B=\{z\in \C~|~y_1<{\rm Im}(z)<y_2\},$$
for some $y_1<y_2$ in $\R$.
\end{LEMMA}

\proof
We will proceed by contradiction. If we cannot put a horizontal band 
in $S$, then there is a horizontal line which intersects both the 
upper boundary of $S$ and the lower boundary of $S$. We can assume 
without loss of generality that this line is the real axis. 
Under the mapping $z\mapsto e^{2i\pi z/\tau}$, the strip projects to 
an annulus $A\in \C^*$, and the real axis projects to the unit circle 
$S^1=\{|z|=1\}$. Hence, the bounded component of $\C^*\setminus A$ contains 
a point $z_1$ of modulus 1, and the unbounded component of 
$\C^*\setminus A$ contains a point $z_2$ of modulus $1$. 

It is known (see \cite{lv} page 56-65) that the modulus of an annulus
separating the points $0$ and $z_1$ from the points $z_2$ and 
$\infty$, with $|z_1|=|z_2|=1$ is bounded from above by the modulus of 
the annulus 
$$A_{max}=\C\setminus\left(]-\infty,-1]\cup [0,1]\right).$$
In particular, when ${\rm mod}(S/\tau)>{\rm mod}(A_{max})$, we see that 
we get a contradiction. 
Douady indicated to us that this modulus is equal to $1/2$. Indeed, 
let us first consider a square pillow with side-length 1. This pillow is 
isomorphic to $\P$. We can map two opposite corners to $0$ and 
$\infty$, and a third corner to $1$.  
By symmetry of the square pillow, the  remaining corner is mapped 
to $-1$. The annulus $A_{max}$ is then isomorphic to the pillow 
cut along two opposite sides. This is a cylinder with height $1$ 
and circumference $2$. Hence the modulus of this cylinder is $1/2$.
\qedlem

\begin{figure}[htbp]
\centerline{
\boxit{\psboxscaled{500}{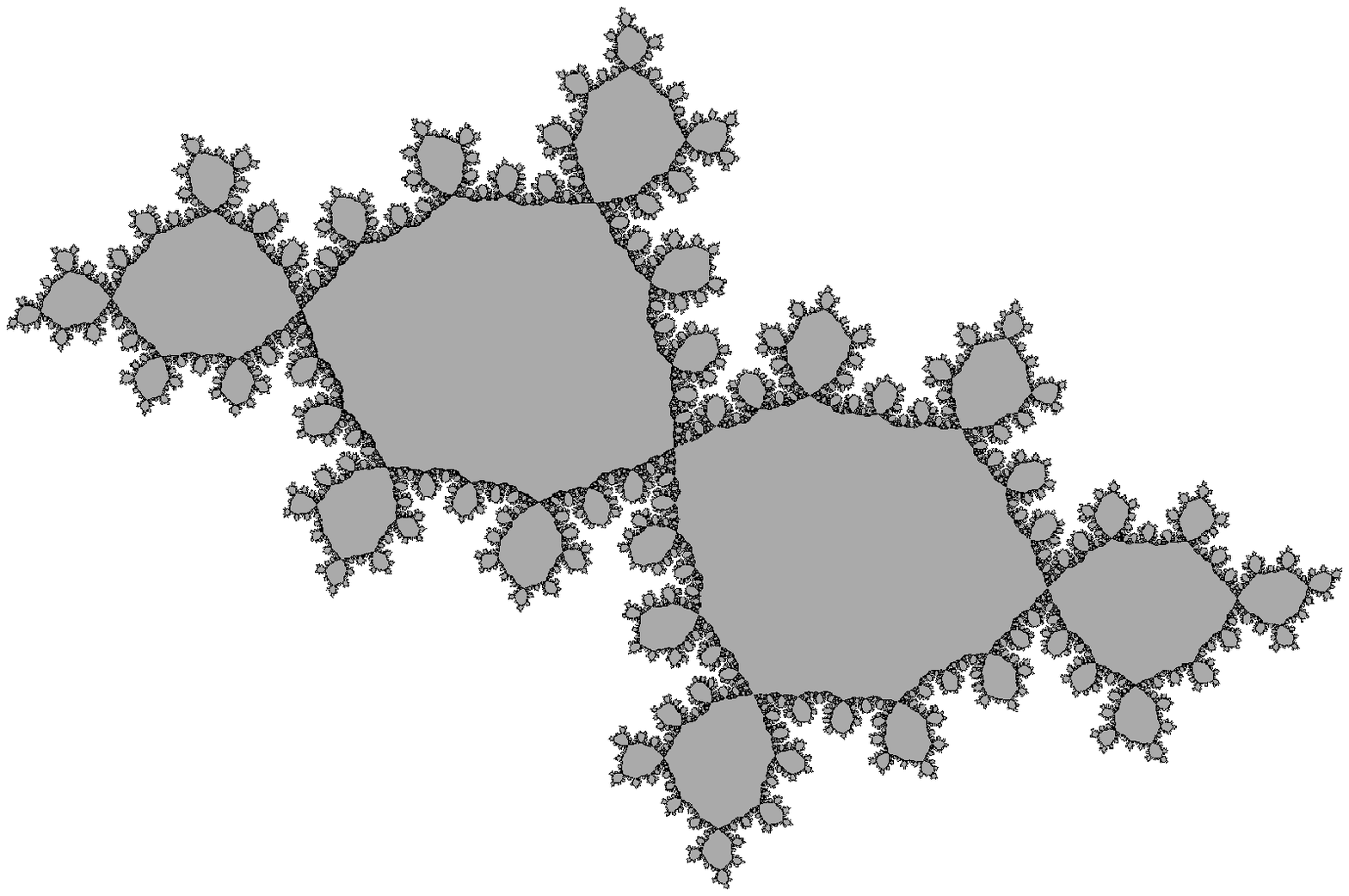}}
\boxit{\psboxscaled{500}{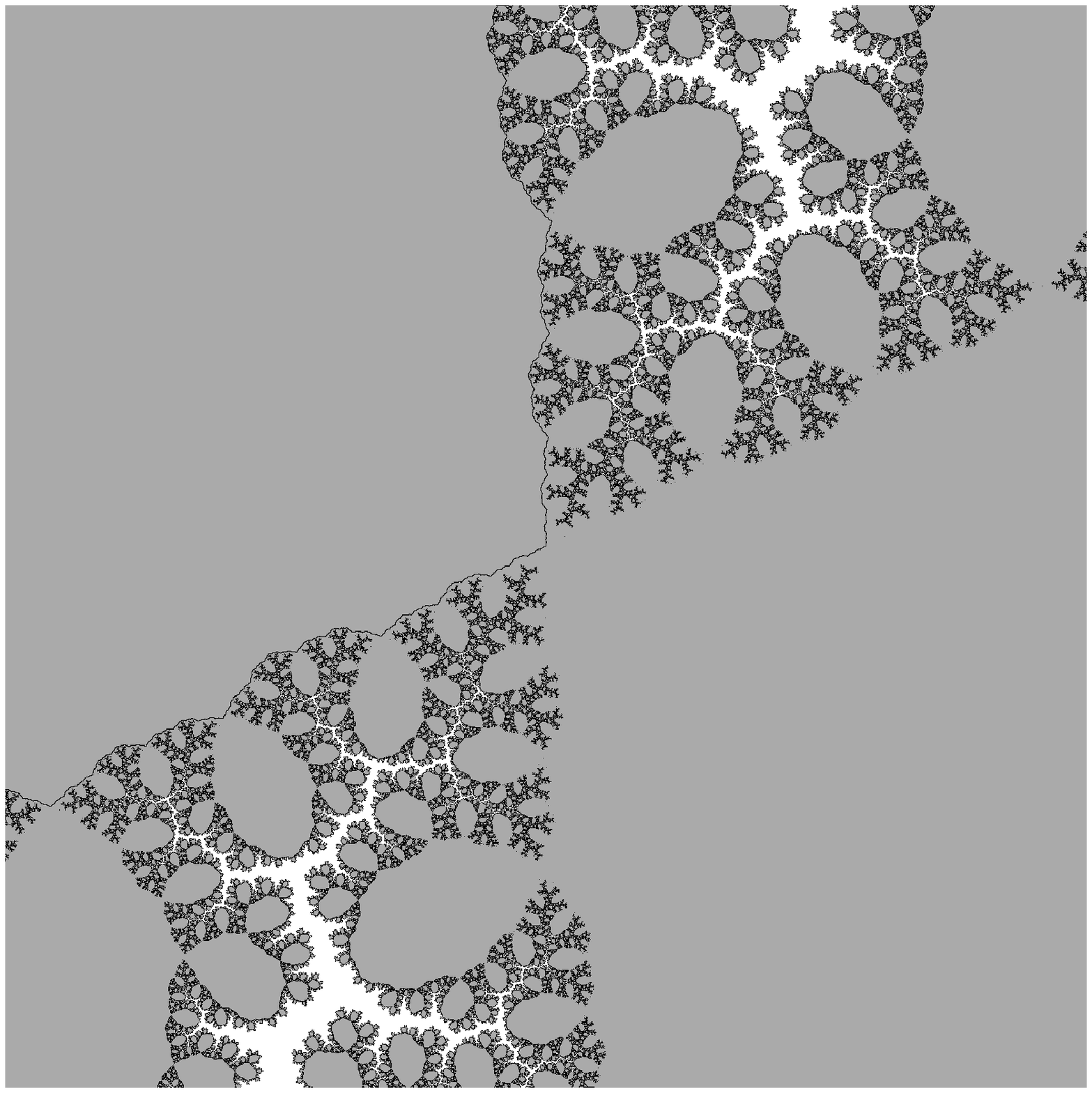}}}
\caption{The filled-in Julia set of the polynomial $P_\theta$, 
$\theta=\sqrt 2-1$.}
\label{figure5}
\end{figure}

\vskip.5cm
\noindent{\bf Remark:} We would like to mention that for the angle
$\theta=[2,2,2,2,\ldots]=\sqrt 2-1$, the modulus of the corresponding 
annulus is 
$$-\pi/\log(\alpha^2)\sim 1.782213977 >\frac{1}{2}.$$
Hence, our proof enables us to conclude
that there is still a triangle in the Siegel disk 
with vertex at the critical point (see figure \ref{figure5}). 

\begin{figure}[htbp]
\centerline{
\boxit{\psboxscaled{500}{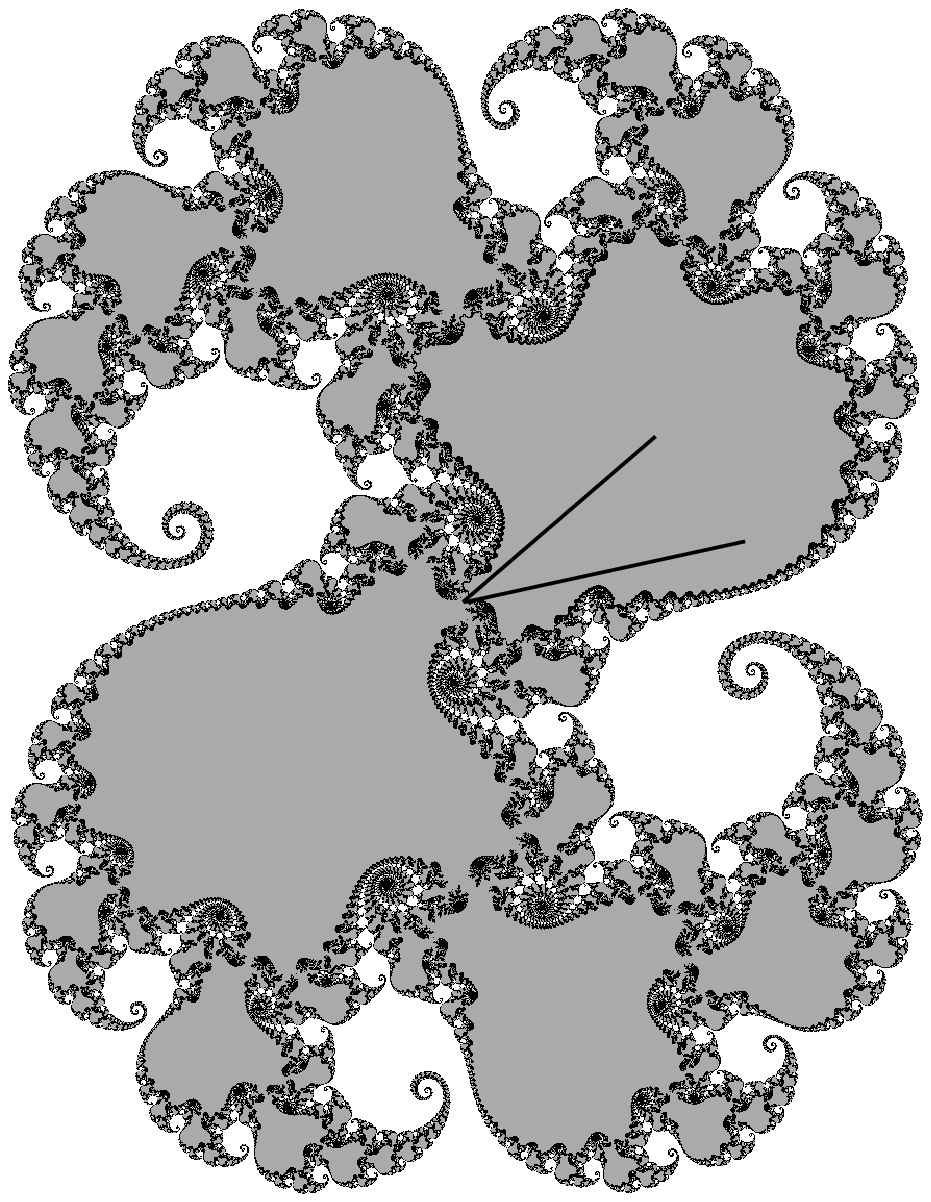}}
\boxit{\psboxscaled{500}{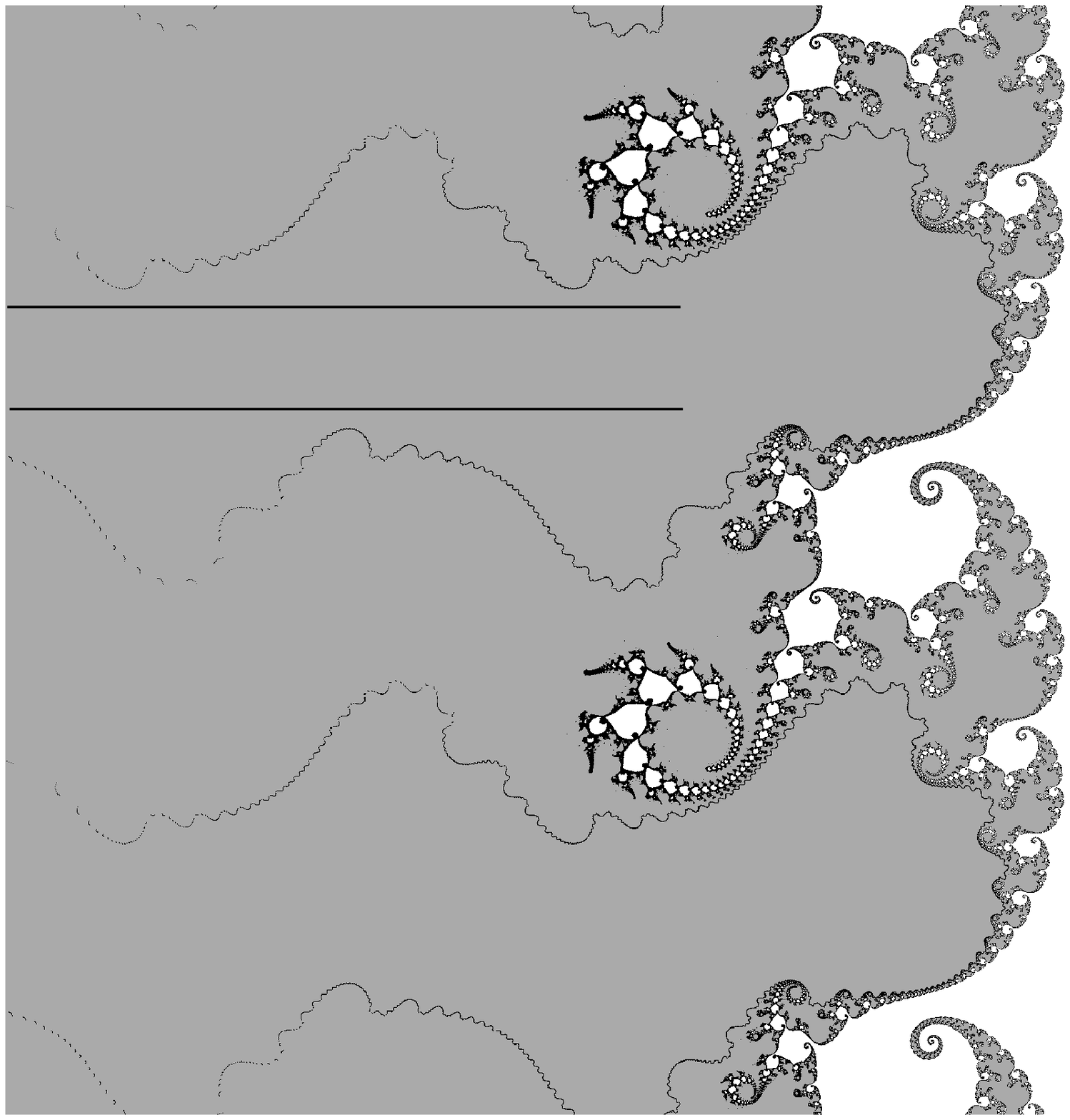}}}
\caption{The filled-in Julia set of the polynomial $P_\theta$, 
$\theta=[24,24,24,\ldots]$ and its image under the map 
$\chi(z)=\log(z-\omega_{\theta})$.}
\label{figure7}
\end{figure}

In fact we can conclude that there is an angle in the Siegel disk as 
long as $\theta=[a,a,a,a,\ldots]$, with $a\leq 23$. Indeed, for 
$a=23$, we have $\alpha=\theta=(\sqrt{533}-23)/2$, and
$$-\pi/\log(\alpha^2)\sim .5006714845 >\frac{1}{2}.$$
On the other hand, when $a=24$, we get $\alpha=\theta=(\sqrt{580}-24)/2$, and
$$-\pi/\log(\alpha^2)\sim .4939944446< \frac{1}{2}.$$
In this case, our proof does not enable us to conclude anything.
We have drawn the corresponding Julia set on figure \ref{figure7}. We 
have also drawn its image under the map 
$\chi(z)=\log(z-\omega_{\theta})$ to show that the boundary of the 
Siegel disk ``oscillates''. It is a reason why our proof does not 
enable us to conclude anything, whereas it seems that one can put a 
triangle in the Siegel disk with vertex at the critical point.

\section{Questions.}

We have seen that when the period of a quadratic irrational is odd, 
then the boundary of the corresponding Siegel disk does not spiral 
around the critical point. 

\vskip.5cm
\noindent{\bf Question 1:} {\em Is there a quadratic irrational 
$\theta$ such that the boundary of the Siegel disk of $P_\theta$ 
spirals?}

\vskip.5cm
\noindent{\bf Question 2:} {\em Does the boundary of the Siegel disk always 
spiral when the period of $\theta$ is even?}

\begin{figure}[htbp]
\centerline{
\boxit{\psboxscaled{500}{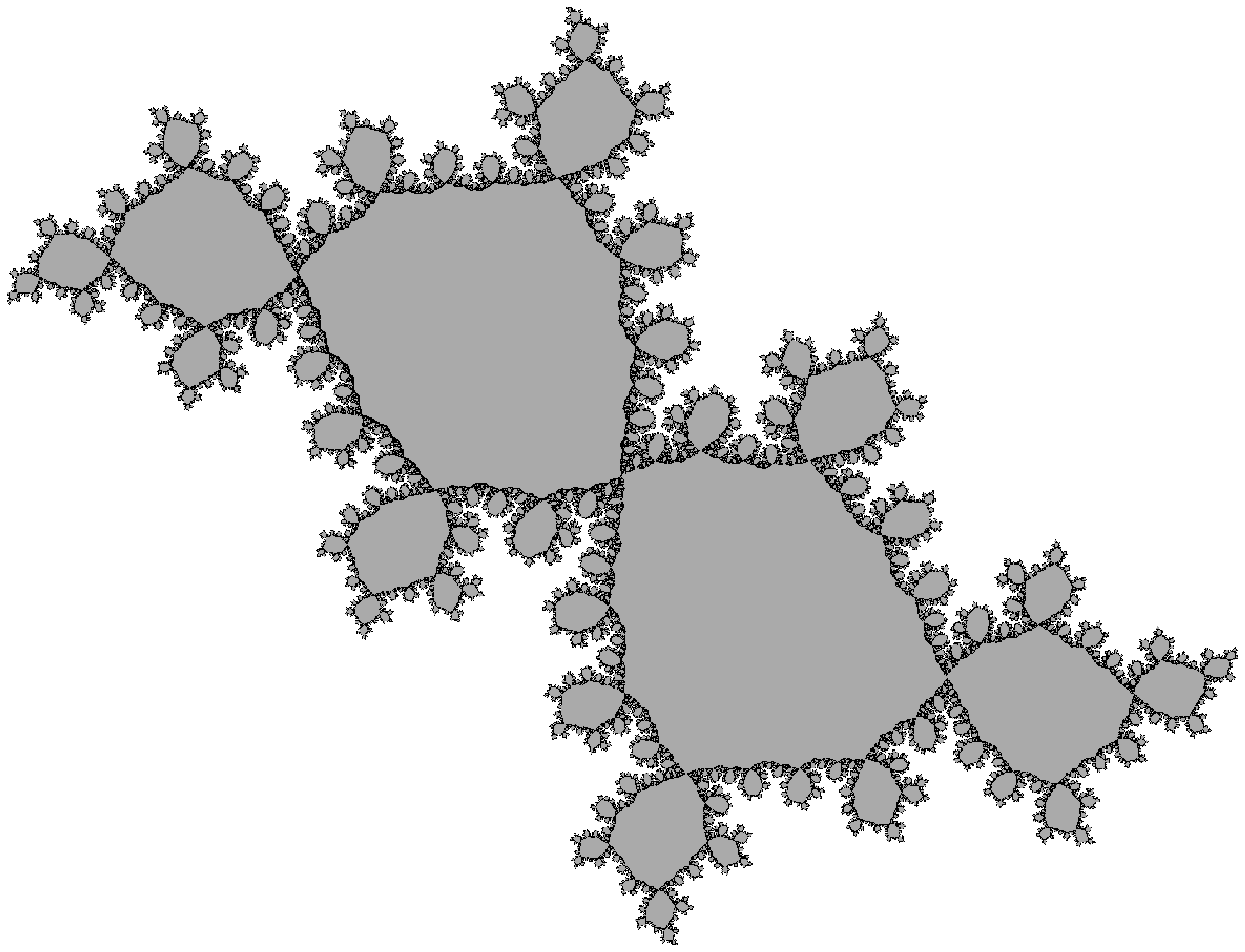}}
\boxit{\psboxscaled{500}{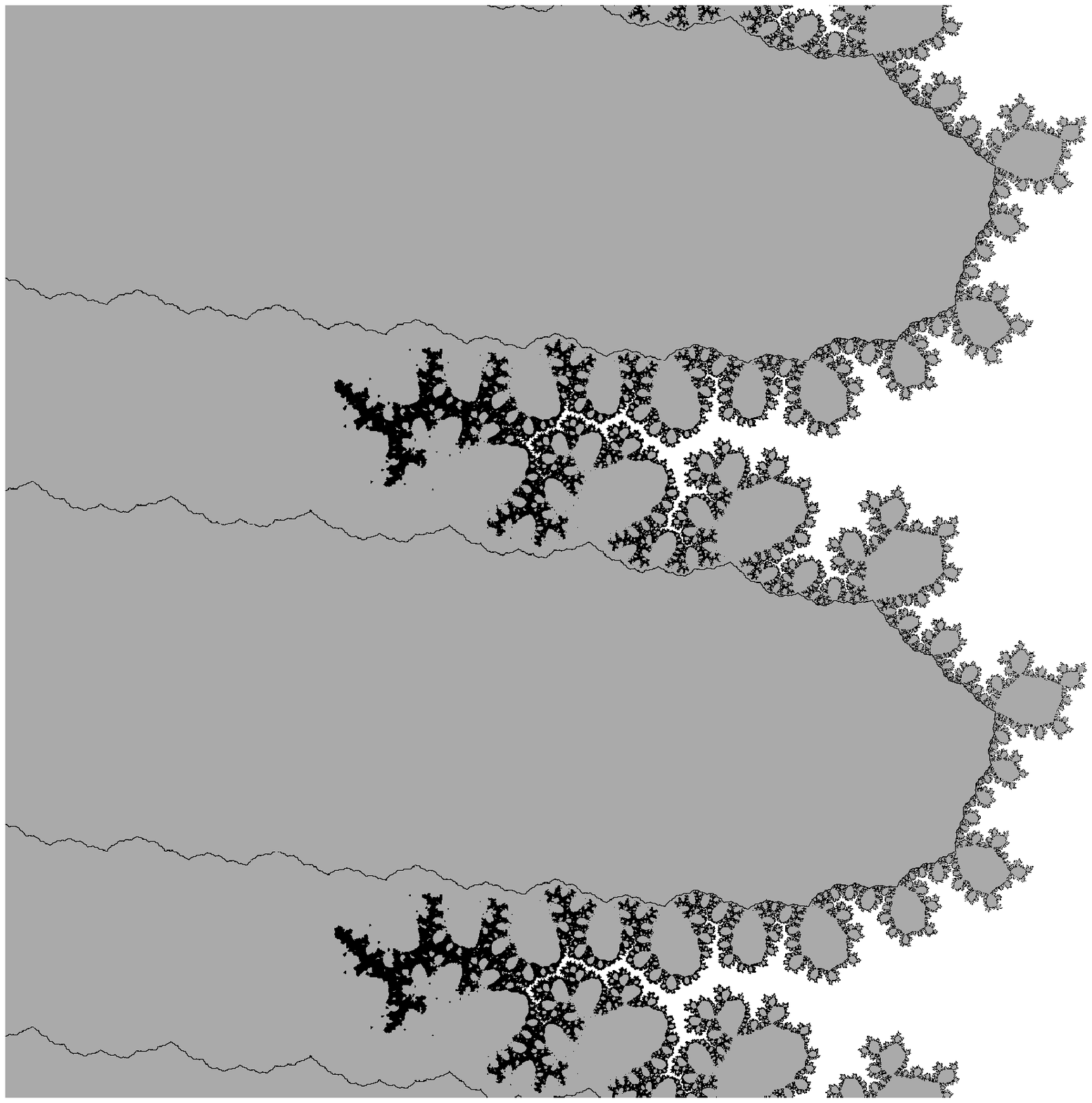}}}
\caption{The filled-in Julia set of the polynomial $P_\theta$, 
$\theta=[2,1,2,1,2,1\ldots]$ and its image under the map 
$\chi(z)=\log(z-\omega_{\theta})$.}
\label{figure9}
\end{figure}

To answer those two questions, one has to show that the scaling 
ratio $\lambda\in\D^*$ is not a real number. 
Computer experiments suggests that for $\theta=[2,1,2,1,2,1,\ldots]$, 
the ratio $\lambda$ is not real. Hence, the boundary of the Siegel 
disk spirals. We have drawn the 
Julia set of the polynomial $P_\theta$ for 
$\theta=[2,1,2,1,2,1,\ldots]$, and its image under the map 
$\chi(z)=\log(z-\omega_{\theta})$ (see figure \ref{figure9}). 
It should be clear that the 
strip corresponding to the Siegel disk is not horizontal.

\vskip.5cm
\noindent{\bf Question 3:} {\em Is there a quadratic irrational 
$\theta$ with odd period, but for which there is no triangle 
with vertex at $\omega_\theta$ contained in the Siegel disk?}

\begin{figure}[htbp]
\centerline{
\boxit{\psboxscaled{500}{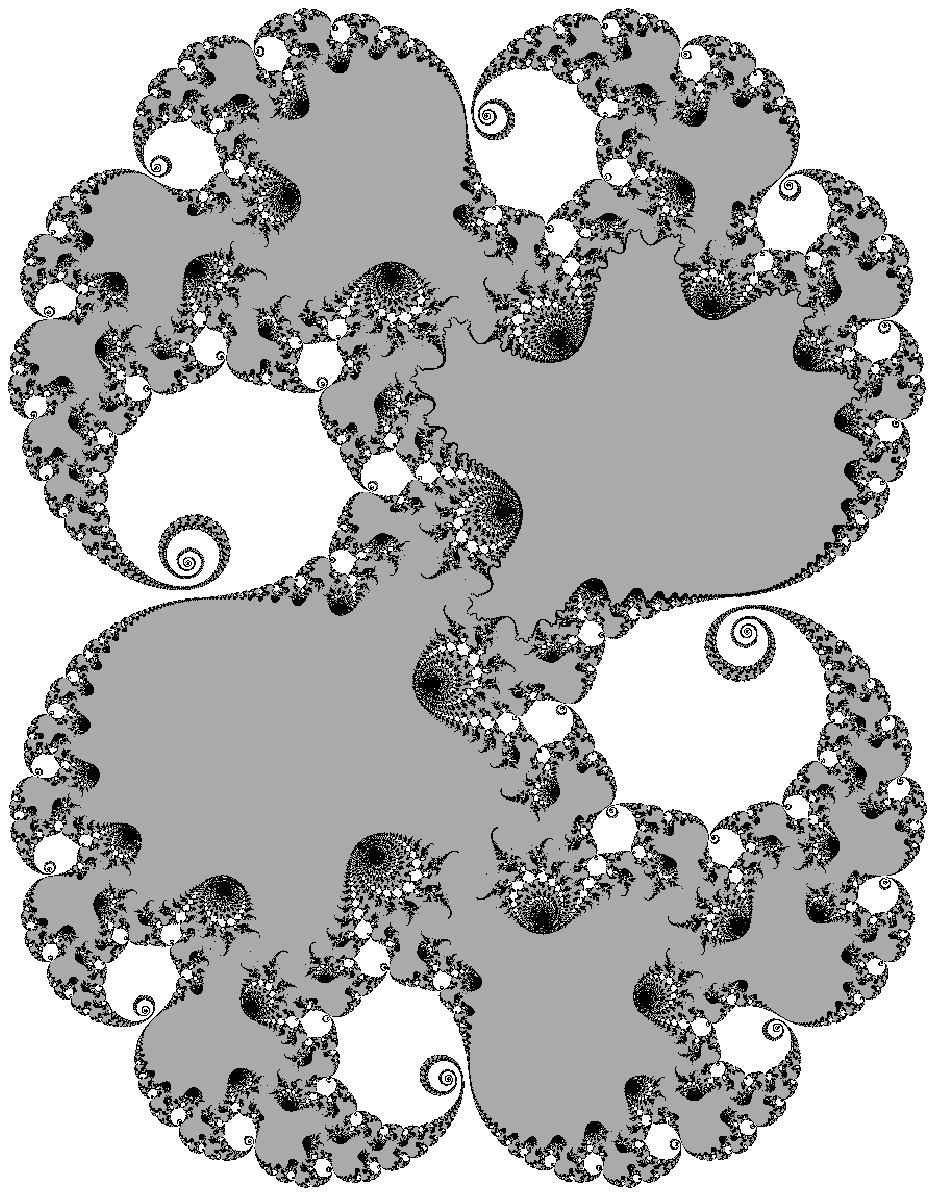}}
\boxit{\psboxscaled{500}{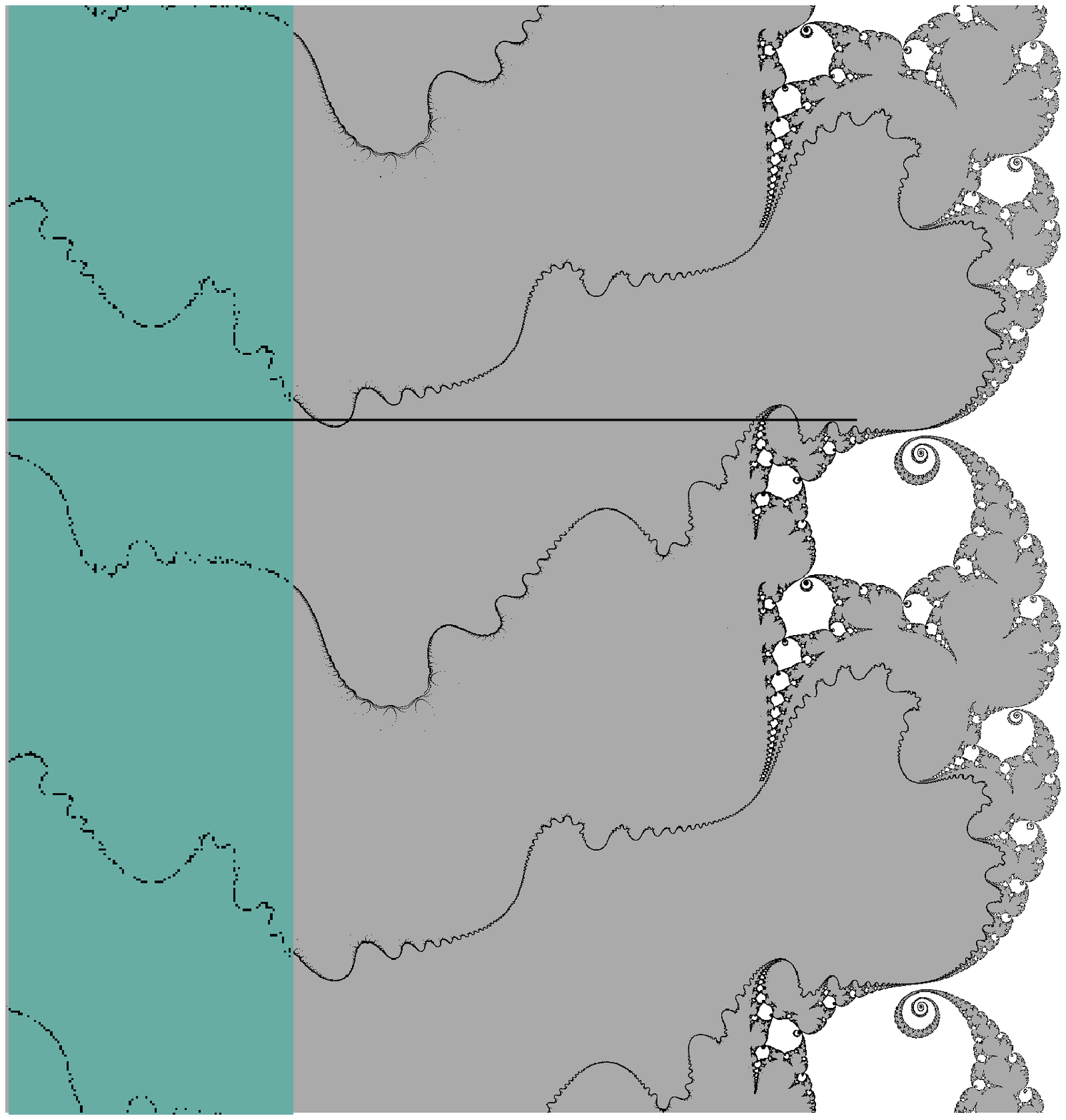}}}
\caption{The filled-in Julia set of the polynomial $P_\theta$, 
$\theta=[50,50,1,50,50,1,50,50,1,\ldots]$ and its image under the map 
$\chi(z)=\log(z-\omega_{\theta})$.}
\label{figure11}
\end{figure}

This problem seems to be related to 
question 1. Indeed, if there is a quadratic irrational such that the 
boundary of the Siegel disk of $P_\theta$ spirals, then the period $s$
of $\theta$ is even. Let us write 
$\theta=[a_1,\ldots,a_s,a_1,\ldots,a_s,\ldots]$. Now consider the 
quadratic irrational $\theta'=[a'_1,\ldots,a'_{ks+1},\ldots]$ 
of period $ks+1$, where $k$ is a large integer, and where
\begin{eqnarray*}
a'_{1} & = & a_1,\quad{\rm and}\\
a'_i & = & a_{i-1},\quad{\rm if}\quad 2\leq i\leq ks+1.
\end{eqnarray*}
Then, one can expect that the boundary of the Siegel disk of 
$P_{\theta'}$ will oscillate. 
On figure \ref{figure11}, we have drawn the filled-in Julia set of 
the polynomial $P_\theta$, 
$\theta=[50,50,1,50,50,1,50,50,1,\ldots]$ and its image under the map 
$\chi(z)=\log(z-\omega_{\theta})$. It is very difficult to obtain a 
good picture of the boundary of the Siegel disk near the critical 
point. The dark region on figure \ref{figure11} corresponds to 
something we extrapolated.

Finally, we know that $|\lambda|<1$ for every quadratic irrational 
$\theta$. When we show that there is a triangle in the Siegel disk, we 
don't get any lower bound on the angle of the vertex which is at the 
critical point. We could get one if we knew that $|\lambda|$ is not 
too close to 1.

\vskip.5cm
\noindent{\bf Question 4:} {\em Is there a constant $\delta<1$ such 
that $|\lambda|<\delta$ for any quadratic irrational $\theta$?}

We can even be more optimistic.

\vskip.5cm
\noindent{\bf Question 5:} {\em Is there a constant $\delta<1$ such 
that $|\lambda|<\delta^s$, where $s$ is the period of 
the quadratic irrational $\theta$?}

Finally, we would like to ask a last question, which seems to be a 
analog of \'Swi\c{a}tek's a-priori bounds for Blaschke fractions 
having an irrational rotation number:

\vskip.5cm
\noindent{\bf Question 6:} {\em Are there constants $\delta_1<\delta_2<1$ such 
that $(\delta_1)^s<|\lambda|<(\delta_2)^s$, where $s$ is the period of 
the quadratic irrational $\theta$?}

\newcounter{nom}{\setcounter{nom}{1}}


\begin{thebibliography}{McM2}

\bibitem[B] {b} {\sc L. Bers}, {\em On Boundaries of Teichm\"uller 
Spaces and on Kleinian Groups}, Annals of Math. {\bf vol} 91 (1970) 570--600.

\bibitem[D] {d} {\sc A. Douady}, {\em Disques de Siegel et Anneaux de Herman},
  S\'eminaire Bourbaki, Ast\'erisque {\bf vol} 152-153 (1986/87) 151--172.

\bibitem[H] {h} {\sc M. Herman}, {\em Conjugaison quasi-sym\'etrique des 
diff\'eomorphismes du cercle et applications aux disques singuliers 
de Siegel}, Manuscript, 1986.

\bibitem[LV] {lv} {\sc O. Lehto} $\&$ {\sc K.I. Virtanen}, {\em 
Quasi-conformal Mappings in the Plane}, Springer-Verlag (1973).

\bibitem[MN] {mn} {\sc N.S. Manton} $\&$ {\sc M. Nauenberg}, 
{\em Universal scaling behavior for iterated maps in the complex plane}, 
Comm. Math. Phys. {\bf 89} (1983), 555--570.

\bibitem[McM1] {mcm1} {\sc C.T. McMullen}, {\em Complex Dynamics and 
renormalization}, Princeton Univ. Press (1994).

\bibitem[McM2] {mcm} {\sc C.T. McMullen}, {\em Self-similarity of
Siegel disks and Hausdorff dimension of Julia sets}, Acta Mathematica,
(1998).

\bibitem[McM3] {mcm3} {\sc C.T. McMullen}, {\em Iteration on Teichm\"uller 
space}, Inv. Math. {\bf vol} 99 (1989), 425--454.
        
\bibitem[P] {p} {\sc C.L. Petersen}, {\em Local connectivity of some
Julia sets containing a circle with an irrational rotation}, Acta 
Math. {\bf 177} (1996), 163--224.

\bibitem[Si] {si} {\sc C.L. Siegel}, {\em Iteration of analytic functions}, 
Annals of math. {\bf 43} (1942), 607--612.

\bibitem[Sw] {sw} {\sc G. \'Swi\c{a}tek}, {\em Rational Rotation Numbers for
    Maps of the Circle}, Comm. Math. Phys. {\bf 119} (1988), 109--128.

\bibitem[W] {w} {\sc M. Widom}, {\em Renormalization group analysis of 
quasi-periodicity in analytic maps}, Comm. Math. Phys. {\bf 
92} (1983), 121--136.

\end{thebibliography}
\end{document}